\magnification=1200
\input amstex
\documentstyle{amsppt}

\pageheight {49pc}
\vcorrection{-2pc}


\def\leftheadline {\line{\hfill\eightrm MARSHALL A. WHITTLESEY\hfill\folio}}
\def\rightheadline{\line{\hfill \eightrm POLYNOMIAL HULLS AND 
H$^\infty$ CONTROL \hfill \folio}}

\headline{\ifodd \pageno\rightheadline \else\leftheadline\fi}
\def\C{\hbox{\bf C$\,$}}
\def\R{\hbox{\bf R$\,$}}
\def\opdisk{\hbox {\rm int$\,\Delta\,\,$}}
\def\Re{\hbox{\rm Re}}
\def\Im{\hbox{\rm Im}}

\topmatter
\title  Polynomial hulls and $H^\infty$ control for 
a hypoconvex constraint
\endtitle
\author Marshall A. Whittlesey
\endauthor
\affil Texas A\&M University
\endaffil
\address College Station, TX 77843-3368
\endaddress
\email mwhittle\@math.ucsd.edu
\endemail
\curraddr Department of Mathematics 0112, University of California, San Diego,
9500 Gilman Drive, La Jolla CA 92093-0112
\endcurraddr
\abstract We say that a subset of $\C^n $ is hypoconvex if its
complement is the union of complex hyperplanes.  Let $\Delta $ be the
closed unit disk in $\C$, $\Gamma=\partial\Delta$.  We prove two
conjectures of Helton and Marshall.  Let $\rho $ be a smooth function
on $\Gamma\times\C^n$ whose sublevel sets have compact hypoconvex
fibers over $\Gamma$.  Then, with some restrictions on $\rho $, if $Y$
is the set where $\rho $ is less than or equal to $1$, the polynomial
convex hull of $Y$ is the union of graphs of analytic vector valued
functions with boundary in $Y$.  Furthermore, we show that the infimum
$\inf_{f\in H^\infty(\Delta)^n}
\|\rho(z,f(z))\|_\infty$ is attained by a unique bounded analytic $f$
which in fact is also smooth on $\Gamma$.  We also prove that if $\rho
$ varies smoothly with respect to a parameter, so does the unique $f$
just found.

\endabstract
\nologo
\endtopmatter

\document
\NoBlackBoxes
\baselineskip=20pt

We address two conjectures of Helton and Marshall from [HMa, p. 183]
which generalize previous theorems regarding an $H^{\infty}$ control
problem over the disk and polynomial hulls of compact sets in
\C$^{n+1}$ fibered over the circle in \C.

If $Y$ is a compact set in $\C^n$, then the {\it polynomial (convex)
hull} $\widehat Y$ of $Y$ is given by $$\hbox {$\widehat Y =
\{z\in\C^n\bigl ||P(z)|\leq
{\displaystyle \sup _{w\in Y}} |P(w)|$ for
all polynomials $P$ on \C$^n\}$.}$$ 
Let $\Delta $ be the closed unit disk in \C , let $\Gamma $ be the
unit circle, let $\Pi:\Delta\times\C\longrightarrow\Delta $ be
projection, and let $\rho:\Gamma\times\C^n\rightarrow \R$ be $C^2$.
Let $D^i\rho=$ the $i^{th}$ derivative of $\rho $, (thought of as an
$i-$linear mapping, as in [L].)  Then $D\rho(z,w) $ and $D^2\rho(z,w)$
are the gradient and real Hessian at $(z,w)$, respectively, for
$(z,w)\in
\Gamma\times\C^n$.  Let $D_w\rho(z,w)$ be the vector $({\partial\rho
\over\partial w_1}(z,w),{\partial\rho
\over\partial w_2}(z,w),...,{\partial\rho
\over\partial w_n}(z,w))$ and let $D_{\overline w}\rho(z,w)$ 
denote its conjugate.
Note that $D_{\overline w}\rho$ is the complex form for the gradient
of $\rho $ in $w$.  In this work, if $L, B$ are linear and bilinear
maps, respectively, we shall write $L[u]$ and $B[u][v]$ to denote the
values of $L$ and $B$.  Where we encounter functions which are not
linear, we shall use parentheses instead of brackets to denote values.
Following [HMa], we define $\rho$ to be {\it hypoconvex} if there exists 
a ${\Cal K}\geq 0$ such that for every
point $(z,w_1,w_2,...,w_n)\in\Gamma\times\C$
such that if $$\sum_{j=1}^n u_j {\partial\rho\over \partial w_j}
(z,w_1,w_2,...,w_n)=0$$ we have
$$D^2\rho(z,w_1,w_2,...,w_n)[0,u_1,u_2,...,u_n][0,u_1,u_2,...,u_n]
\geq \Cal K\|u\|^2,\tag 1$$
where $\|\cdot\|$ is the standard Euclidean norm.  Let
$K^t=\{(z,w_1,w_2,...,w_n)\in\Gamma\times\C^n\bigl|\rho(z,w_1,w_2,$
$...,w_n)=t\}$ and let $K^t_z=\{(w_1,w_2,...,w_n)\bigl|
(z,w_1,w_2,...,w_n)\in K^t\}$ be the fiber of $K^t$ over $z$.  If
${\Cal K}>0$ then (1) says that on the complex tangent space of the
fibers, $D^2\rho $ is positive definite, so the complex hyperplane in
\C$^n$ tangent to the set $K^t_z$ is locally external to the set where
$\rho$ is less than or equal to $t$.

In [H\" o], a set in $\C^n$ is defined to be {\it linearly convex} if
the complement of the set is the union of complex $(n-1)$-dimensional
affine planes.  The notion of such sets appeared in 1935 under the
name ``planarkonvex'' in work of Behnke and Peschl [BP] (although the
notion is slightly weaker), again in 1940 in [BS], later discussed in
the 1952 dissertation [St] of Strehlke, and then reappeared in the
1960's in work of Martineau, who used the term ``lin\' eellement
convexe'' (see [M]).  Kiselman (see [Ki]) uses the term ``lineal
convexity'' and Vityaev (see [V]), ``complex geometric convexity.''
In order to avoid confusion with ordinary convexity, we shall use the
word hypoconvex instead of linearly convex.  We also feel the word
hypoconvex is suggestive of the geometry of the situation, since this
notion is somewhat weaker than ordinary convexity.  However, the
reader should be aware of the fact that we are not following the
majority of the literature in using this terminology.

In summary, we shall call a set in $\C^n$ {\it hypoconvex} if its
complement is the union of complex ($n-1$)-dimensional affine planes
and we shall call $\rho $ hypoconvex on an open subset $U$ of
$\Gamma\times\C$ if it satisfies (1) in $U$. If for every compact
subset of $U$ there exists a ${\Cal K}>0$ such that (1) holds, we
shall say that $\rho $ is {\it strictly hypoconvex} on $U$.  A
$C^2$-bounded hypoconvex set whose defining function has real Hessian
which is positive definite on the complex tangent spaces at every
boundary point of the set shall be called {\it strictly hypoconvex}.
Thus $K^t_z$ bounds a strictly hypoconvex compact set for all
$z\in\Gamma $.

We shall work under the assumption that $n\geq 2$ but our
arguments work for $n=1$ with minor adjustments.

Our first plan is to prove (in Theorem 2) that a compact set $K$ in
$\Gamma\times\C^n$ with smoothly bounded connected strictly hypoconvex
fibers containing the origin has a polynomial hull which is the union
of the graphs of analytic vector-valued functions over the closed unit
disk, provided that the set can be deformed in a reasonable manner to
a compact set whose fibers are balls.  This has been proven provided
$K$ has convex fibers in [AW] and [S3], and if $n=1$ for compact $K$
with connected and simply connected fibers in [F1] and [S2], and also
in [HMa].  Our result then generalizes both, since convex sets in
$\C^n$ are hypoconvex, and any subset of \C is hypoconvex.  There
exist examples of compact sets in
\C$^n$ fibered over the circle with contractible fibers whose
polynomial hulls are not the union of graphs over the disk; see [HMe1]
and [\v Ce].  We say that the function $f$ defined on $\Gamma $ is a
{\it selector} for the set $K\subset\Gamma\times\C^n$ if $f(z)\in K_z$
for all $z\in\Gamma $.

A problem of $H^\infty$ control is to compute
$$\gamma_\rho\equiv\inf_{f\in H^\infty(\Delta)^n}\hbox{\rm
ess}\sup_{z\in\Gamma}
\rho(z,f(z)),\tag 2$$
using notation from [HMa].  We shall call $\gamma_\rho $ the {\it
optimal control} for $\rho $.  It is also of interest to determine
various facts about an $f$ which attains this minimax: whether it
exists, is unique, is smooth, and whether it possesses other
properties to be mentioned later.  Such an $f$ we call a ``solution''
to the $H^\infty$ control problem (2).  We shall prove (in Theorem 3)
that for strictly hypoconvex $\rho $ the $H^\infty$ control problem is
uniquely solvable with a vector valued function which is smooth, again
subject to restrictions to be described.  This has been done in some
important cases in [HMa],[Hu] and [HV].  [HV] uses methods similar to
the ones we employ and proves related results.  In [V], Vityaev shows
that if $\rho $ is strictly hypoconvex where $\rho=\gamma_\rho$, then
if a solution to (2) is smooth on $\Gamma $, it is the only smooth
solution.

We also show in Theorem 4 that if $\rho $ changes smoothly with
respect to a parameter, then so does $\gamma_\rho$ and the solution
to (2).

We shall make the following assumptions on $\rho$ throughout our work:
$$\eqalign {(3)&\left [\vbox {\hsize =4.5 in\vbox{(a)
$\rho:\Gamma\times\C^n\rightarrow[0,\infty)$ is continuous, and
$C^6$-smooth where $\rho\neq 0$;}\medskip
\vbox {(b) there exists an $R>1$ such that if $0<\rho(z,w)\leq R$
then $D_w\rho(z,w)\neq 0$ and $\rho$ is strictly hypoconvex as 
in (1) where $\rho $ is smooth;}\medskip
\vbox{(c) for every $t$, $0<t\leq R$, the set $K^t$ where $\rho=t$ is
compact, with fibers $K^t_z$ diffeomorphic to $\partial
B_n$, where $B_n$ is the open unit ball in $\C^n$;}
\vbox{(d) $K=K^1$;}
\vbox{(e) $\{(z,w) \bigl| \rho(z,w)=R\}=\{(z,w) \bigl | |w|=R\}$ and 
$\rho(z,w)>R$ if $|w|>R$;}\medskip
\vbox{(f) There exists a continuous function $S(z)$, $|S(z)|<R$, 
such that $\{(z,w)\in
\Gamma\times\C^n\bigl| w=S(z)\}=\{(z,w)\in
\Gamma\times\C^n\bigl| \rho(z,w)=0\}$.}\vskip -1.25in}\right .}$$

From Theorem 1 of [YK], we may conclude that if $0<t<R$, and $D^t_z$
is the closed domain enclosed by the level set
$\{w\,|\,\rho(z,w)=t\}$, then the complex tangent spaces to boundary
points of $D^t_z$ do not meet $\hbox{\rm int}\,D^t_z$.  See also [Ki].
By Corollary 4.6.9 of [H\" o] applied to $\hbox{\rm int}\, D^t_z$, we
may conclude that $\hbox{\rm int}\,D^t_z$ is hypoconvex, so $D^t_z$ is
as well, as the intersection of $D^s_z$ for $s>t$.  $D^t_z$ is also
clearly strictly hypoconvex.  From the same
corollary, we obtain that $\hbox{\rm int}\,D^t_z$ is ``\C-convex'',
i.e., if $P$ is a 1-dimensional complex affine subspace of $\C^n$
which meets $\hbox{\rm int}\,D^t_z$, then the intersection is a
connected and simply connected subset of $P$.  If $P\cap\hbox{\rm
int}\,D^t_z\neq\emptyset$, then the intersection is a smoothly bounded
subset of $P$, because $P$ cannot be a tangent to $D^t_z$; hence the
derivative of $\rho $ restricted to $P$ where $\rho=t$ is not
identically zero and $P\cap D^t_z$ is the closure of $P\cap\hbox{\rm
int}\,D^t_z$.  If  $P\cap\hbox{\rm int}\,D^t_z=\emptyset$ but $P$
does meet $D^t_z$ then we claim it only meets $D^t_z$ in one point.
Such a $P$ must be tangent to the boundary of $D^t_z$ at any intersection
point.  As observed earlier, the tangent space is locally disjoint from
$D^t_z$ near such a point.  If there are two such points, then
we may perturb $P$ slightly and obtain a $P'$ whose intersection
with int $D^t_z$ is not connected.  Thus the intersection of any complex
hyperplane with $D^t_z$ is either a point or a Jordan domain with
boundary as smooth as $\rho$.  We note that $D^t_z\subset
\widehat{K^t_z}$.  The reverse inclusion also holds: Proposition 1 of
[Z] shows that the interior of $D^s_z$ is polynomially convex for all
$s$, so $D^t_z$ is polynomially convex as the intersection of the
polynomially convex open sets $D^s_z$, $s>t$.  Thus $D^t_z=\widehat
{K^t_z}$.

\S 1 {\bf A perturbation theorem.} 

Let $\langle\cdot,\cdot\rangle $ be the complex inner product on
\C$^n$, $\langle a,b\rangle =\sum_{j=1}^n a_j\overline {b_j}$.  We
shall allow the arguments $a$ and $b$ to be functions whose values are
in \C$^n$; then the operation $\langle\cdot,\cdot\rangle $ is
pointwise inner product of the two functions.  Let
$A(\Delta)=\{f:\Delta\rightarrow\C |
\hbox{$f$ is continuous on $\Delta$ and analytic in \opdisk}\}$,
$L^2_C(\Gamma)=$ the set of square integrable complex valued functions
on $\Gamma$ under the ordinary inner product $\langle f,g\rangle
_2=\int_0^{2\pi} f(e^{i\theta})\overline
{g(e^{i\theta})}\,{d\theta\over 2\pi}$, $W^{1,2}(\Gamma)=$ the Sobolev
space of complex functions on $\Gamma $ with first derivatives in
$\theta$ in $L^2_C(\Gamma)$, under the real inner product $\langle
f,g\rangle _{1,2}=C\Re\langle f,g\rangle _2+\Re\langle {\partial
f\over
\partial\theta},{\partial g\over \partial\theta}\rangle _2,$ where $C$
is a large positive constant to be chosen later.  Let
$W^{1,2}_R(\Gamma)=$ the real valued elements in $W^{1,2}(\Gamma)$,
$H^2(\Delta)=$ the elements of $L^2_C(\Gamma)$ whose negative Fourier
coefficients vanish and $H^{1,2}(\Delta )=H^2(\Delta)\cap
W^{1,2}(\Gamma)$.  Then $W^{1,2}(\Gamma)\subset C(\Gamma)$ is a
continuous inclusion by the Sobolev embedding theorem; but we can show
quickly why this holds.  If $f(z)=\sum_{j=-\infty}^\infty a_jz^j$ is
in $W^{1,2}(\Delta )$ then $\|f_\theta\|_2^2=\sum_{j=-\infty}^\infty
j^2|a_j|^2<\infty$ so $\sum_{j=-\infty}^\infty |a_j|\leq
|a_0|+(\sum_{j=-\infty}^\infty j^2|a_j|^2)^{1\over 2}
(2\sum_{j=1}^\infty{1\over j^2})^{1\over 2}=|a_0|+{\pi\over\sqrt{3}}
\|f_\theta\|_2$.  Thus the Fourier series of $f$ converges absolutely
and uniformly on $\Gamma $ and its supremum norm is bounded by $\sum
_{j=-\infty}^\infty|a_j|\leq |a_0|+{\pi\over\sqrt{3}}\|f_\theta\|_2$,
which is an equivalent norm for $W^{1,2}$.  Then $f\in A(\Delta )$ and
has small supremum norm if $\|f\|_{1,2}$ is small.  We let
$H^{1,2}_0(\Delta )=$ the set of elements of $H^{1,2}(\Delta )$ which
have value zero when $z=0$.  

In the following we use techniques similar to those of Forstneri\v c
[F1].  We begin with the graph of an analytic vector valued
function $f$ which is extremal in a sense that will be clear
later (it will be in the boundary of a particular polynomial
hull) and find graphs close by which are similarly extremal.
We spell out conditions we shall require of such graphs
and use the implicit function theorem on Banach spaces to 
establish their existence and uniqueness.

{\bf Theorem 1.}  {\it Let $\rho $ satisfy (3) with $S$ identically
zero.  Suppose that there exist $f,g\in H^{1,2}(\Delta)^n$ such that
$f,g\in C^4(\Gamma)$, $\rho(z,f(z))$ is constant in $z$, say $=1$,
${\displaystyle\sum_{j=1}^{n}}f_j(z)g_j(z)=1$ and the affine complex
tangent plane to $K^1_z$ at $(z,f(z))$ is $\{(w_1,w_2,...,w_n)\in\C
^n\bigl| \sum_{j=1}^n g_j(z)w_j=1\}$.  Then for some neighborhood
$N(f(0))$ in \C$^n$, there exist $C^1$ maps $F,G:N(f(0))
\rightarrow H^{1,2} (\Delta)^n$, $F=(F_1,F_2,...,F_n)$, 
$G=(G_1,G_2,...,G_n)$, such that $\rho(z,F(w)(z))$ is 
constant in $z\in\Gamma$ for fixed $w\in N(f(0))$, $F(w)(0)=w$,
$\sum_{j=1}^{n}F_j(w)(z)G_j(w)(z)=1$, and the complex tangent plane to
$K^{\rho(z,F(w)(z))}_z$ at $F(w)(z)$ is $\{(w_1,w_2,...,w_n)\in\C
^n\bigl|\sum_{j=1}^n G_j(w)(z)w_j=1\}$.}

{\it Remark.}  For Theorem 1, it suffices to assume that $\rho $ is
only $C^4$ and that the set where $\rho $ equals $R$ bounds
a set which is merely strictly convex instead of being a ball.

{\it Proof.}  Condition (3)(f) with $S=0$ guarantees that $K^t_z$
separates the origin in $\C^n$ from the point at infinity for $t>0$.

We first reduce to the case where $g(z)=(1,0,0,...,0)$ and $f_1$, the
first coordinate of $f$, is identically 1.  To do this, we construct
an $n\times n$ matrix $M(z)$ of analytic functions such that the first
column is given by $g(z)$.  By [SW], Theorem 2.1, since $g(z)$ is
never zero on $\Delta $, there exist analytic $\C^n$- valued functions
$h_1$, $h_2$,...,$h_{n-1}$ in $A(\Delta )^n$ such that $g(z)$ along
with the $h_i$ generate $A(\Delta )^n$ as a module over $A(\Delta )$.
In particular, if the last $n-1$ columns of $M$ are
$h_1,h_2,...,h_{n-1}$, then for all $z\in\Delta $, the determinant of
$M(z)$ is nonzero.  We may mollify the $h_i$ so slightly they are also
in $C^4(\Gamma )$ but the determinant of $M(z)$ is still nonzero.
Then consider the $C^4$ function given by $\tilde\rho(z,w)
=\rho(z,w\cdot M(z)^{-1})$, regarding $w$ as a row vector.  Then we
get corresponding sets $\tilde K^t\equiv \{(z,w)\in\Gamma
\times\C^n|\tilde\rho(z,w)=t\}$ such that under $w\mapsto w\cdot M(z)$,
the complex tangent space to $K_z$ at $f(z)$ is mapped to the complex
tangent space to $\tilde K_z$ at $(1,0,0,...,0)$, which is
$\{(w_1,w_2,...,w_n)|w_1=1\}$.  Let $\tilde g(z)= (1,0,0,...,0)$.  We
also note that we get an $\tilde f(z)\equiv f(z)\cdot M(z)$
corresponding to $f(z)$ such that $\tilde f_1(z)=1$.  Then the
conditions for the theorem are satisfied with $\tilde\rho$, $\tilde f$
and $\tilde g$, except that the conditions of the Remark above hold.
The linearity and invertibility of $w\mapsto w\cdot M(z)$ in $w$
guarantees that $\tilde\rho $ is strictly hypoconvex.  Once we have
Theorem 1 proven for $\tilde\rho$, obtaining $\tilde F(w),\tilde G(w),
\tilde N(\tilde f(0))$ with the desired properties, we define
$F(w)(z)=\tilde F(w\cdot M(0))(z)\cdot M(z)^{-1}$.  Then $F$ is
defined in a neighborhood $N(f(0))$ of $f(0)$ since $\tilde
f(0)=f(0)\cdot M(0)$.  We have $\rho(z,F(w)(z))=\rho(z,\tilde F(w\cdot
M(0))(z)\cdot M(z)^{-1})=\tilde\rho(z,\tilde F(w\cdot M(0))(z))$ which
is constant in $z$ for fixed $w$ in some neighborhood of $f(0)$.  If
we write $v=w\cdot M(z)$, then ${\partial\rho\over\partial w_j}=$
$\sum_{i=1}^n ({\partial\tilde\rho\over\partial v_i}{\partial
v_i\over\partial w_j}$ $+ {\partial\tilde\rho\over\partial
\overline{v_i}}{\partial
\overline{v_i}\over\partial w_j})$= $\sum_{i=1}^n
{\partial\tilde\rho\over\partial v_i} {\partial v_i\over\partial
w_j}$.  Let $x^T$ denote the transpose of $x$.  Now ${\partial
v_i\over\partial w_j} =M_{ji}$, so
$D_w\rho(z,F(w)(z))^T=M(z)(D_w\tilde\rho(z,\tilde F(w\cdot
M(0))(z)))^T= a(w)(z)M(z)\cdot (\tilde G(w\cdot M(0))(z))^T$, where
$a(w)$ is $C^1$ in $w$ and we claim the winding number of $a(w)$ is
zero in $z$.  If $w=f(0)$ then this is merely the statement that
$\sum_{j=1}^n\tilde f_j(z){\partial\tilde\rho\over\partial
w_j}(z,\tilde f(z))$ has winding number zero.  This function is not
zero for $z\in\Gamma$ because the complex tangent planes to the fibers
of $\tilde K^t$ never pass through the origin.  Now we may deform
$\tilde f$ through a homotopy$\{\tilde f^t\}$, $1\leq t\leq R$ so that
$\tilde f^t$ is continuous, $\tilde f^1=\tilde f$, $\rho(z,\tilde
f^t(z))=t$ and $\sum_{j=1}^n\tilde f^t_j(z){\partial\tilde\rho
\over\partial w_j}(z,\tilde f^t(z))$ has the same winding number as
$\sum_{j=1}^n\tilde f_j(z){\partial\tilde\rho\over\partial
w_j}(z,\tilde f(z))$; the former is never zero for $z\in\Gamma $ for
the same reason as stated above for $\sum_{j=1}^n\tilde
f_j(z){\partial\tilde\rho\over\partial w_j}(z,\tilde f(z))$.  Then
$\tilde K^R$ can be deformed smoothly out to some $\tilde K^{R'}$
which has spherical fibers of radius $\tilde R'$ such that $\tilde
K^t$ has convex fibers for $R\leq t\leq R'$.  Deforming $\tilde f^R$
similarly out to $\tilde f^{R'}$ such that $\rho(z,\tilde f^t(z))=t$
for $R\leq t\leq R'$, we find that $\sum_{j=1}^n\tilde
f^{R'}_j(z){\partial\tilde\rho\over\partial w_j}(z,\tilde
f^{R'}(z))=|\tilde f^{R'}(z)|^2={R'}^2,$ a constant function (with
winding number $0$).  By the homotopy to $\sum_{j=1}^n\tilde f_j(z)
{\partial\rho\over\partial w_j}(z,\tilde f(z))$, $a(v)(z)$ has winding
number $0$ for $v=f(0)$ and for all $v$ in $N(f(0))$.
Then $D_w\rho(z,F(v)(z))$ is a scalar function $a(v)(z)$ times an
analytic vector function $b(v)(z)$ which is $C^1$ in $v$, so
$\sum_{j=1}^nF_j(v)(z)b_j(v)(z)$ is analytic and equals ${1\over
a(v)(z)}$ times $\sum_{j=1}^nF_j(v)(z){\partial\rho\over \partial
w_j}(z,F(v)(z))$, both of which have winding number zero.  (The latter
function has winding number zero by an argument in $K^t$ similar to
the one just given in $\tilde K^t$.)  We can then define
$$G(w)(z)\equiv {b(w)(z)\over \sum_{j=1}^nF_j(w)(z)b_j(w)(z)}$$ for
$w$ in $N(f(0))$.  Then $F(w)$ and $G(w)$ satisfy the properties
required.  This proves Theorem 1 provided we prove it in the case
where $g(z)=(1,0,0,...,0)$ and $f_1=1$.

We note that there exists a function
$h\in H^{1,2} (\Delta)$ such that $h$ is
nonzero for $z\in\Delta$ and $\sum_{j=0}^nf(z){\partial \rho\over
\partial w_n}(z,f(z))={\partial\rho\over\partial w_1}(z,f(z))$ 
has the same argument as $h(z)$ on $\Gamma$.  To obtain $h$,
we need the fact that this function is nonzero on $\Gamma$, with
winding number $0$, using an argument from the previous paragraph.  We
also need that the harmonic conjugation operator is continuous on
$W^{1,2}$, that if $f\in W^{1,2}$ then so is $e^f$, and that $\log(f)$
is as well if $f$ is nonzero and has winding number $0$.

Suppose that $u,v\in W^{1,2}_R(\Gamma).$ Let $\tilde u$ denote the
harmonic conjugate of $u$ whose value at $0$ is $0$.  Let $H^{1,2}
(\Delta)^{n-1}$ be the subspace of $H^{1,2} (\Delta)^n$ of $n$-tuples
$k=(k_1,k_2,...,k_n)$ such that $k_1=0$.  By $W^{1,2}_R(\Gamma)/\R$ we
mean the quotient of $W_R^{1,2}$ by the real constant functions.

As in [F1], let $X(z)$ denote an element of $H^{1,2}(\Delta )$ which
points in the same real direction as the outward normal to $K^1_z$ at
$1$, which is ${\partial\rho\over\partial\overline{w_1}}(z,f(z)).$ (In
fact, we can take $X=1/h$.)  Then wind $X(z)=0$, so $X$ has no zeroes
on the closed disk.  Let
$$F(u,v,k)=f+(u+\tilde
ui)fX+v(0)ifX+k$$
$$G(l)=g+l$$
and $F=(F_1,F_2,...,F_n)$, $G=(G_1,G_2,...,G_n)$.  Also let $w_0=f(0).$

Consider the mapping 
$\Phi:W^{1,2}_R(\Gamma)\times W^{1,2}_R(\Gamma)\times H^{1,2}
(\Delta)^{n-1}\times H^{1,2}(\Delta)^n\times \C^n
\rightarrow W^{1,2}_R(\Gamma)/\R\times W^{1,2} (\Gamma)^n\times
\C^n$, where

\noindent
$\Phi(u,v,k,l,w)=(\Phi_1(u,v,k,l,w),\Phi_2(u,v,k,l,w) -
\Phi_3(u,v,k,l,w),\Phi_4(u,v,k,l,w))$
and

\noindent
$\Phi_1(u,v,k,l,w)(\cdot)=\rho(\cdot,F(u,v,k)(\cdot))+\R,$

\noindent
$\displaystyle {\Phi_2(u,v,k,l,w)(z)={1+v(z)+i\tilde v(z)\over
\sum_{j=1}^n (\overline {F_j(u,v,k)(z)}) ({\partial\rho\over\partial
\overline {w_j}}(z,F(u,v,k)(z)))}}\bullet$

$\hskip 0.5in\biggl(\displaystyle {{\partial \rho\over\partial 
\overline{w_1}}
(z,F(u,v,k)(z)),{\partial \rho\over\partial
\overline{w_2}}(z,F(u,v,k)(z)),...,{\partial \rho\over\partial
\overline {w_n}}(z,F(u,v,k)(z))\biggr),}$

\noindent
$\Phi_3(u,v,k,l,w)(z)=\bigl(\overline {G_1(l)}(z),
\overline {G_2(l)}(z),...,\overline {G_n(l)}(z)\bigr),$ and

\noindent
$\Phi_4(u,v,k,l,w)=\int_0^{2\pi}F(u,v,k)(e^{i\theta}){d\theta\over 2\pi}-w.$

Note that $\Phi(0,0,0,0,w_0)=(0,0,0)$.  In \S 5, we show that $\Phi$
is a $C^1$ map near $(0,0,0,0,w_0)$.  We claim that the partial
derivative of $\Phi$ in $(u,v,k,l)$ when $(u,v,k,l,w)=(0,0,0,0,w_0)$
is an invertible mapping.  We denote this partial derivative by
$D_{(u,v,k,l)}$.  Using the implicit function theorem, we will be able
to make the conclusions of the theorem.  The reader is invited to look
at the end of the proof to see how this happens.  We also note
that our technique resembles that used by Lempert in [L2].

Following Forstneri\v c [F1], we compute
$D\Phi_1(0,0,0,0,w_0)[u,0,0,0,0]=$ \linebreak $2\Re
(\sum_{j=1}^n{\partial\rho\over\partial w_j}(\cdot,f(\cdot))(u+\tilde
ui)f_jX)+\R=2\Re({\partial\rho\over\partial
w_1}(\cdot,f(\cdot))(u+\tilde ui)X)+\R$.  Since $X(z)$ has the same
argument as ${\partial\rho\over\partial{\overline w_1}}(z,f(z))$,
$2\Re({\partial\rho\over\partial w_1}(\cdot,f(\cdot))(u+\tilde
ui)(X)=$\linebreak $2{\partial\rho\over\partial
w_1}(\cdot,f(\cdot))X\Re(u+\tilde ui)=2{\partial\rho\over\partial
w_1}(\cdot,f(\cdot))Xu$.  Then $I\equiv 2{\partial\rho\over\partial
w_1}(\cdot,f(\cdot))X$ is a positive $W_R^{1,2}$ function.

We proceed with several steps, initially showing that
$D_{(u,v,k,l)}\Phi(0,0,0,0,w_0)$ is injective in $(u,v,k,l)$.  First
we note that $D\Phi_1(0,0,0,0,w_0)[0,v,k,l,0]=0$ for all $v\in
W^{1,2}_R(\Gamma)$, $k\in (H^{1,2} (\Delta))^{n-1}$ and
$l\in (H^{1,2} (\Delta))^n$.  This is obvious with $l$.  For $k$ and
$v$ we see that
$D\Phi_1(0,0,0,0,w_0)[0,v,k,0,0](z)=D\rho(z,f(z))[v(0)if(z)X(z)+k(z)]=0,$
since $v(0)if(z)X(z)+k(z)$ is a tangent to $K_z$ at $(z,f(z))$.

We shall have need for the fact that $\langle
D\Phi_2(0,0,0,0,w_0)[0,0,k,0,0],f\rangle =0$.  To see this, we first
note that $\langle \Phi_2(u,v,k,l,w_0),F(u,v,k)\rangle =1+v+\tilde
vi$.  Differentiating both sides with respect to $k$ at $(u,v,k,l,w)=
(0,0,0,0,w_0)$, we find that
$$\langle D\Phi_2(0,0,0,0,w_0)[0,0,k,0,0],F(0,0,0)\rangle +\langle
\Phi_2(0,0,0,0,w_0), DF(0,0,0)[0,0,k]\rangle =0.$$ Now
$DF(0,0,0)[0,0,k]=k$, $F(0,0,0)=f$, $\Phi_2(0,0,0,0,w_0)=(1+v+\tilde
vi)
\overline g$ and $\langle \overline g,k\rangle =0$ so the above
simplifies to
$$\langle D\Phi_2(0,0,0,0,w_0)[0,0,k,0,0],f\rangle =0,\tag 4$$
as desired.

Clearly $D\Phi_2(0,0,0,0,w_0)[0,0,0,l,0]=0$ and
$D\Phi_3(0,0,0,0,w_0)[0,0,k,0,0]=0$ for all $k\in H^{1,2}
(\Delta)^{n-1}$ and $l\in H^{1,2} (\Delta)^n$, as $\Phi_2$ is
constant in $l$ and $\Phi_3$ is constant in $k$.

Suppose that $D\Phi(0,0,0,0,w_0)[u,v,k,l,0]=0$.  Then $uI$ is a
constant $c$, where $I$ is a nonzero real $W^{1,2}$ function.  In
order for $D\Phi_4(0,0,0,0,w_0)[u,v,k,l,0]$ to equal zero, we must
have $0=\Re \langle D\Phi_4(0,0,0,0,w_0)[u,v,k,l,0],X(0)\overline
{g(0)}\rangle = u(0)|X(0)|^2$.  Since $X$ has winding number zero, we
have $u(0)=0.$ Hence $c\bigl(({1\over I})(0)\bigr)=0$, where (${1\over
I})(0)$ denotes the value at $0$ of the harmonic extension of $1/I$ to
the closed disk.  Now $({1\over I})(0)\neq 0$ since $I$ is positive,
so we must have $c=0$.  Thus $u=0$.  From the fact that
$D\Phi_4(0,0,0,0,w_0)[0,v,k,l,0]=0$, we also find that $v(0)=0$
($0=\langle D\Phi_4(0,0,0,0,w_0)[0,v,k,l,0],X(0)\overline g(0)\rangle
=v(0)|X(0)|^2i$) and
$k(0)=\int_0^{2\pi}k(e^{i\theta})\,d\theta={\buildrel \rightarrow\over
0}$, where ${\buildrel \rightarrow\over 0}=$ the origin in $\C^n$.

Then the derivative of the middle coordinate of $\Phi $ must be $0$ in
the direction of $[0,v,k,l,0]$, so
$$\eqalign {0=&D\Phi_2(0,0,0,0,w_0)[0,v,k,l,0]-D\Phi_3(0,0,0,0,w_0)
[0,v,k,l,0]\cr =&D\Phi_2(0,0,0,0,w_0)[0,v,k,0,0]-D\Phi_3(0,0,0,0,w_0)
[0,0,0,l,0],\cr}$$ making use of the previously observed facts that
$D\Phi_2(0,0,0,0,w_0)$ does not depend on $l$ and
$D\Phi_3(0,0,0,0,w_0)$ does not depend on $v,k$.  Thus
$$\langle D\Phi_2(0,0,0,0,w_0)[0,v,k,0,0],f\rangle =
\langle D\Phi_3(0,0,0,0,w_0)[0,0,0,l,0],f\rangle .$$ 
From (4), we get $$\langle
D\Phi_2(0,0,0,0,w_0)[0,v,0,0,0],f\rangle=\langle
D\Phi_3(0,0,0,0,w_0)[0,0,0,l,0],f\rangle ,$$ so then $1+v+\tilde
vi=(1+v+\tilde vi)\langle \overline g,f\rangle =\langle
D\Phi_2(0,0,0,0,w_0)[0,v,0,0,0],f\rangle$ $ =\langle
D\Phi_3(0,0,0,0,w_0)[0,0,0,l,0],f\rangle =\sum_{j=1}^n\overline
{l_i}\overline {f_i}$.  (Note that we need the fact that $v(0)=0$
implies $DF(0,0,0)[0,v,0]=0$.)  Thus $v+\tilde vi=0$ (as an analytic
and conjugate analytic function whose value at $0$ is $0$.)

So $\langle D\Phi_2(0,0,0,0,w_0)[0,0,k,0,0],k\rangle =$ $\langle
D\Phi_3(0,0,0,0,w_0)[0,0,0,l,0],k\rangle $.  Now suppose $k$ is not
identically $0$; then
$$\eqalign {\langle &D\Phi_2 (0,0,0,0,w_0)[0,0,k,0,0],k\rangle (z)\cr=& 
\langle {D^2\rho(z,f(z))[k(z)] \over\sum_{j=1}^n (\overline
{f_j}(z))({\partial\rho\over\partial
\overline {w_j}}(z,f(z)))},k(z)\rangle \cr&+C(z)\sum_{j=1}^n
{\partial\rho\over
\partial \overline {w_1}}(z,f(z))\overline{k_j}(z)\cr
&\hbox {(using the fact that $v=0$)}\cr
=&\langle {D^2\rho(z,f(z))[k(z)]
\over{\partial\rho\over\partial
\overline {w_1}}(z,f(z))},k(z)\rangle \cr}$$
since $k(z)$ is a complex tangent to $K_z$ at $f(z)$.  However we also
have \linebreak$\langle D\Phi_3(0,0,0,0,w_0)[0,0,0,l,0],k\rangle (z)=$
$\sum_{j=1}^n
\overline {k_j(z)}\overline {l_j(z)}$ so 
$$\overline {h(z)}\langle {D^2\rho(z,f(z))[k(z)]
\over{\partial\rho\over\partial
\overline {w_1}}(z,f(z))},k(z)\rangle =\sum_{j=1}^n
\overline{h(z)}\overline {k_j(z)}\overline {l_j(z)}$$
so, taking real parts of both sides and integrating,
$$\int_0^{2\pi}\hskip -2 mm\Re{\overline
{h(z)}\over {\partial\rho\over\partial
\overline {w_1}}(z,f(z))}D^2\rho(z,f(z))[k(z)][k(z)]\,
{d\theta\over 2\pi}=\int_0^{2\pi}\hskip -2 mm\Re\sum_{j=1}^n
\overline{h(z)}\overline {k_j(z)}\overline {l_j(z)}
\,{d\theta\over 2\pi}.\tag 5$$
The right side of (5) is $0$ since $k_j(0)=0$ for all $j$.  However
the integrand on the left side is positive since $${\overline
{h(z)}\over {\partial\rho\over\partial\overline {w_1}}(z,f(z))}$$ is
positive, by definition of $h$, and since the real Hessian of $\rho $
in $w$ is positive definite on complex tangents on the fibers of $K$.
Thus $k=0$, so $D\Phi_3(0,0,0,0,w_0)[0,0,0,l,0]$ $=0$, hence $l=0$.

In a similar manner, we can show that $D_{(u,v,k,l)}\Phi(0,0,0,0,w_0)$
is a surjective linear map.  Let $(u'+\R,m,c)$ be an element in the
target space.  Then $D\Phi(0,0,0,0,w_0)[u'/I,0,0,0,0]$ has first
coordinate $u'+\R$, so we may assume without loss of generality that
$u'=0$ in $W_R^{1,2}/\R$, but in the domain we must restrict ourselves
to the subspace where $u=b/I$ for some real constant $b$.  Similarly,
if 
$$b={\Re\, \langle c,X(0)\overline{g(0)}\rangle \over ({1\over
I}(0))|X(0)|^2}$$ then $\Re\langle
D\Phi_4(0,0,0,0,w_0)[b/I,0,0,0,0],X(0)\overline {g(0)}\rangle =\Re\,
\langle c,X(0)\overline{g(0)}\rangle $.  Hence we may assume that
$\Re\, \langle c,X(0)\overline{g(0)}\rangle =0$, provided that we
restrict the domain to the subspace in the domain where $b=0$ (i.e.,
$u=0$).  Next, if we let $v=$ the constant function $\Im\,\langle
c,X(0)\overline{g(0)}\rangle /(|X(0)|^2)$, then $$\Im\langle
D\Phi_4(0,0,0,0,w_0)[0,v,0,0,0],X(0)\overline {g(0)}\rangle
=\Im\langle c,X(0)\overline{g(0)}\rangle $$ so we may assume $0=\langle
c,X(0)\overline{g(0)}\rangle =\overline{X(0)}c_1$, so $c_1=0$,
provided we restrict the domain to the subspace where $v(0)=0$.  Then
we must show that $D\Phi_4(0,0,0,0,w_0)$ maps $H^{1,2}(\Delta)^{n-1}$
onto the subspace of $\C^n$ of complex dimension $n-1$ of $n$-tuples
with first coordinate zero; this is obvious.  Thus we now restrict our
domain to those $k\in H^{1,2}(\Delta )^{n-1}$ for which
$k(0)={\buildrel \rightarrow\over 0}$.

We must then show that the image of $$\{0\}\times\{v\in
W^{1,2}_R(\Delta)|v(0)=0\}\times
\{k\in H^{1,2}(\Delta)^{n-1}|k(0)={\buildrel \rightarrow\over 0}\}\times 
H^{1,2}(\Delta)^n\tag 6$$ under
$D_{(u,v,k,l)}\Phi_2(0,0,0,0,w_0)-D_{(u,v,k,l)}\Phi_3(0,0,0,0,w_0)$,
is onto $W^{1,2} (\Delta)^n$.  First we show the image dense and
then show it to be closed.  Suppose that $m\in W^{1,2}
(\Delta)^n$ is orthogonal to the image.  Then clearly $m\in H^{1,2}
_0(\Delta)^n$.  Further, $m$ is orthogonal to $(1+v+\tilde
vi,0,0,...,0)$ for all $v\in W^{1,2}_R(\Delta)$ with $v(0)=0$, so
orthogonal to all $(p,0,0,...,0)$ where $p\in H^{1,2}_0(\Delta)$; thus
$0=m_1,$ the first coordinate of $m$.  Then in $W^{1,2} (\Delta)^n$,
$m$ is real orthogonal to $D\Phi_2(0,0,0,0,w_0)[0,0,k,0,0]$ for all
$k\in H^{1,2}_0(\Delta)^{n-1}$ with $k(0)={\buildrel \rightarrow\over
0}$.  Let $k=mX.$ Then
$$\eqalign {0=&\langle D\Phi_2(0,0,0,0,w_0)[0,0,mX,0,0],m\rangle
_{1,2}\cr =&C\Re\langle D\Phi_2(0,0,0,0,w_0)[0,0,mX,0,0],m\rangle
_2\cr &+\Re\langle {\partial\over\partial \theta}
(D\Phi_2(0,0,0,0,w_0)[0,0,mX,0,0],{\partial\over\partial
\theta}m\rangle _2.
\cr}\tag 7$$
Now $D\Phi_2(0,0,0,0,w_0)$ acts pointwise on elements in $H^{1,2}_0(
\Delta)^{n-1}$; we may represent its action by $$\eqalign{k(z)\mapsto 
C(k)(z)(1,0,0,...,0)&+{1\over {\partial\rho\over\partial
\overline {w_1}}(z,f(z))}\biggl(\tilde k(z)\Psi(z)\biggr),}$$ 
where $\Psi(z)$ is
the real Hessian of $\rho $ and $\tilde k(z)=(\Re k_1(z),\Im
k_1(z),\Re k_2(z),\Im k_2(z),...,$ $\Re k_n(z),\Im k_n(z))$ is
the realification of $k(z)$.  We must complexify $\tilde k(z)\Psi(z)$
before multiplication by ${1\over {\partial\rho\over\partial\overline
{w_1}}(z,f(z))}$.  Thus (7) equals
$$\eqalign{&C\Re\langle  {1\over {\partial\rho\over\partial
\overline{w_1}}(\cdot,f(\cdot))}   \widetilde{mX}\Psi,m\rangle _2\cr &+
\Re\langle {\partial\over\partial
\theta}  \left({1\over {\partial\rho\over\partial
\overline{w_1}}(\cdot,f(\cdot))} (\widetilde{mX}\Psi)\right),
{\partial\over\partial \theta}m\rangle _2\cr &\hbox{(using the fact
that the first coordinate of $m$ is zero)}\cr=&C\Re\langle
\widetilde{mX}\Psi,{mX\over {\partial\rho\over\partial
w_1}(\cdot,f(\cdot))X}\rangle _2\cr &+\Re\langle
{\partial\over\partial
\theta}  \left({1\over {\partial\rho\over\partial
\overline{w_1}}(\cdot,f(\cdot))} (\widetilde{mX}\Psi)\right),
{\partial\over\partial \theta}(mXX^{-1})\rangle _2\cr=& C\Re\langle
\widetilde{mX}\Psi,{mX\over {\partial\rho\over\partial
w_1}(\cdot,f(\cdot))X}\rangle _2\cr &+\Re\langle
{\partial\over\partial\theta}({1\over {\partial\rho\over\partial
\overline{w_1}}(\cdot,f(\cdot))}) (\widetilde{mX}\Psi)+({1\over 
{\partial\rho\over\partial
\overline{w_1}}(\cdot,f(\cdot))})  ((\widetilde{mX})_\theta\Psi)\cr
&+
({1\over {\partial\rho\over\partial
\overline{w_1}}(\cdot,f(\cdot))}) (\widetilde{mX}\Psi_\theta),{\partial\over
\partial\theta}(mX)X^{-1}+(mX) {\partial\over
\partial\theta}X^{-1}\rangle _2\cr}$$
$$\eqalign{
=& C\Re\langle \widetilde{mX}\Psi,{mX\over {\partial\rho\over\partial
w_1}(\cdot,f(\cdot))X}\rangle _2 \cr &+\Re\langle ({1\over
{\partial\rho\over\partial
\overline{w_1}}(\cdot,f(\cdot))}  ((\widetilde{mX})_\theta\Psi), 
{\partial\over
\partial\theta}(mX)X^{-1}\rangle _2\cr
&+\hbox{\rm other terms}\cr 
=& C\Re\langle \widetilde{mX}\Psi,{mX\over {\partial\rho\over\partial
w_1}(\cdot,f(\cdot))X}\rangle _2 \cr  
&+ \Re\langle   ((\widetilde{mX})_\theta\Psi),{(mX)_\theta\over {\partial
\rho\over\partial
w_1}(\cdot,f(\cdot))X}\rangle _2\cr&+\hbox{\rm other terms.}\cr}\tag
8$$ 
We recall that ${\partial\rho\over\partial w_1}(\cdot,f(\cdot))X$
must be a real-valued function.  Since $\Psi (z)$ is positive definite
on the complex tangent space to $K_z$ at $(z,f(z))$, the first two
terms in this last sum are greater than or equal to ${\Cal
K'}(C\|mX\|_2^2+\|(mX)_\theta\|_2^2)
\geq {\Cal K''}(C\|m\|_2^2+\|m_\theta\|_2^2)$.  
(Recall that the real inner product of $a$ and $b$ is $\Re \langle
a,b\rangle $.  See also the argument leading up to (5).)  The last
group of terms can be written in modulus as less than or equal to
$C_2\|m\|_2\|m_\theta\|_2\leq
C_3\|m\|_2^2/\epsilon^2+C_3\|m_\theta\|_2^2\epsilon^2$, where $C_3$
doesn't depend on $m$.  Choose $\epsilon $ so small that
$C_3\epsilon^2<{1\over 2}{\Cal K}''$ and then assume that $C$ was
chosen large enough that the last equality of (8) is greater than or
equal to a constant times ${\Cal K}/2(C_1\|m\|_2^2
+\|m_\theta\|_2^2)\geq C'\|m\|^2_{1,2}.$ We may then conclude that
$m=0$.  Hence in $(W^{1,2}(\Gamma))^n$, the image of (6) under
$D_{(u,v,k,l)}\Phi_2(0,0,0,0,w_0)-D_{(u,v,k,l)}\Phi_3(0,0,0,0,w_0)$ is
dense.

We need to show that the image of (6) under
$D_{(u,v,k,l)}(\Phi_2-\Phi_3)(0,0,0,0,w_0)$ is closed.  Consider a
convergent $\{m^j\}\in W^{1,2}(\Delta )^n$ in its image; then there
exist $\{v^j\}$, $\{k^j\}$ and $\{l^j\}$ with $k^j(0)={\buildrel
\rightarrow\over 0}$ and $v^j(0)=0$ such that $$m^j=
D\Phi(0,0,0,0,w_0)[0,v^j,k^j,l^j,0].$$  Then
$$m^a-m^b=(1+v^a-v^b+\tilde{v^a}i-\tilde{v^b}i)\overline
g+D\Phi_2(0,0,0,0,w_0) [0,0,k^a-k^b,0,0]+\overline l^a-\overline l^b$$
so, using (4),
$$\langle m^a-m^b,f\rangle
=1+v^a-v^b+\tilde{v^a}i-\tilde{v^b}i+\langle \overline l^a-
\overline l^b,f\rangle $$
converges in $W^{1,2}(\Delta )$.  Projecting the right side to $H^{1,2}$,
we conclude that $\{v^a+\tilde {v^a}i\}$ converges, so $\{v^a\}$
does as well in $W^{1,2}(\Delta )$.  Thus we may assume without
loss of generality that $v^j=0$.
Thus 
$$m^a-m^b=D\Phi_2(0,0,0,0,w_0) [0,0,k^a-k^b,0,0]+\overline l^a-
\overline l^b,$$ 
so $$\langle m^a-m^b,(k^a-k^b)/X\rangle _{1,2}=\langle
D\Phi_2(0,0,0,0,w_0) [0,0,k^a-k^b,0,0], (k^a-k^b)/X\rangle _{1,2}.$$
By reasoning above in (8), the right hand side is $\geq
C_4\|k^a-k^b\|^2_{1,2}$ and the left hand side is $\leq
C_5\|m^a-m^b\|_{1,2}\|k^a-k^b\|_{1,2}$.  Hence $\|k^a-k^b\|_{1,2}\leq
C_6\|m^a-m^b\|_{1,2}$.  We conclude that $\{k^a\}$ converges in
$H^{1,2}(\Delta )$ to some $k$, and similarly that $\{l^a\}$ does.
This proves closure of the image.

From the conclusions of the last several paragraphs, we conclude that
the partial derivative $D_{(u,v,k,l)}\Phi(0,0,0,0,w_0)$ is surjective.

We conclude by the implicit function theorem that there exist an
$N(w_0)$ and unique mappings $u,v,k,l:N(w_0)\rightarrow \bigl(
W^{1,2}_R (\Gamma)\times W^{1,2}_R(\Delta) \times H^{1,2}
(\Delta)^{n-1}\times H^{1,2}(\Delta))^n\bigr)$ such that for
$w\in N(w_0)$, $$\Phi_1(u(w),v(w),k(w),l(w),w)=0+\R,$$
$$\Phi_3(u(w),v(w),k(w),l(w),w)-(1+v(w)+\widetilde{
{v(w)}}i)\Phi_2(u(w),v(w),k(w),l(w),w)=0,$$
and $$\Phi_4(u(w),v(w) ,k(w),l(w),w)=w$$ so $(1+v(w)+\widetilde
{v(w)}i)\Phi_2(u(w),v(w),k(w),l(w),w)=\Phi_3(u(w),v(w),k(w),l(w),w)$.
We also find functions $F(w)\equiv F(u(w),v(w),k(w))$ and $\tilde
G(w)\equiv G(l(w))$ in the same small neighborhood $N(w_0)$, finding
that $F(w)(0)=w$; for fixed $w\in N(w_0)$,
$$\rho(z,F(w)(z))$$ is constant
and for fixed $w$
$${(1+v(w)(z)+\widetilde {v(w)}(z)i)\biggl(\displaystyle{\partial
\rho\over\partial
\overline{w_1}}
(z,F(w)(z)),{\partial \rho\over\partial
\overline{w_2}}(z,F(w)(z)),...,{\partial \rho\over\partial
\overline {w_n}}(z,F(w)(z))\biggr)\over
\displaystyle\sum_{j=1}^n (\overline {F_j}(w)(z)) ({\partial\rho\over\partial
\overline {w_j}}(z,F(w)(z)))}$$ $$=\bigl(\overline {\tilde G_1(w)}(z),
\overline {\tilde G_2(w)}(z),...,\overline {\tilde G_n(w)}(z)\bigr),$$
so
$${(1+v(w)(z)-\widetilde {v(w)}(z)i)\biggl(\displaystyle{\partial
\rho\over\partial {w_1}} (z,F(w)(z)),{\partial \rho\over\partial
{w_2}}(z,F(w)(z)),...,{\partial \rho\over\partial
{w_n}}(z,F(w)(z))\biggr)\over
\displaystyle\sum_{j=1}^n ({F_j}(w)(z)) ({\partial\rho\over\partial
{w_j}}(z,F(w)(z)))}$$ $$=\bigl({\tilde G_1(w)}(z), {\tilde
G_2(w)}(z),...,{\tilde G_n(w)}(z)\bigr),$$ which means that the
complex tangent space to $K^{\rho(z,F(w)(z))} _z$ at $F(w)(z)$ is
$$\{(w_1,w_2,...,w_n)\bigl| \sum_{j=1}^n \tilde G_j(w)(z)w_j=1+v(w)-
\widetilde {v(w)}i\}$$
and
$$\sum_{j=1}^n F_j(w)\tilde G_j(w)=1+v(w)-\widetilde {v(w)}i.$$ Since
the left hand side is analytic and the right hand side is conjugate
analytic, the right side is constant (for fixed $w$); this constant is
$1$ for $w=w_0$ so is nonzero for $w$ near $w_0$.  Letting
$G(w)=\tilde G(w)/(1+v(w)-\widetilde {v(w)}i)$ for $w$ near $w_0$, we
have $F,G$ satisfying the requirements of the theorem.  This concludes
the proof of Theorem 1.$\square$

\S 2 {\bf Extension of the implicit functions $F$ and $G$.}

For $t>0$, we let $L^t$ be the compact set which is the image of $K^t$
under the mapping
$$\eqalign {\Gamma\times\C^n&\longrightarrow\Gamma\times\C^{2n}\cr
(z,w_1,w_2,...,w_n)&\longmapsto (z,w_1,w_2,...,w_n,I_1(z,w),I_2(z,w),
...,I_n(z,w)),}$$
where $w=(w_1,w_2,...,w_n)$ and
$(I_1(z,w),I_2(z,w),...,I_n(z,w))=$ $$\left({{\partial\rho\over
\partial w_1}(z,w)\over \sum_{j=1}^nw_j {\partial\rho\over
\partial w_j}(z,w)},{{\partial\rho\over
\partial w_2}(z,w)\over \sum_{j=1}^nw_j {\partial\rho\over
\partial w_j}(z,w)},...,{{\partial\rho\over
\partial w_n}(z,w)\over \sum_{j=1}^nw_j {\partial\rho\over
\partial w_j}(z,w)}\right).$$
Then $L^t$ is a $C^5-$smooth manifold 
embedded in \C$^{2n+1}$, as it is a $C^5$-smooth graph over $K^t$.
We shall show that $L^t$ is totally real and that the accumulation
points over $\Gamma$ of the analytic disks parametrized by 
$$\eqalign {\opdisk&\longrightarrow \C^{2n+1}\cr z&\longmapsto 
(z,F(w)(z),G(w)(z))}$$
lie in $L^t$.  Results of \v Cirka [\v Ci] will then show that $F(w)$
and $G(w)$ are in $C^4(\Gamma)$ for $w\in N(w_0)$.

The following lemma is closely related to a lemma in [Wb].

{\bf Lemma 1.}  {\it $L^t$ is totally real.}

{\it Proof.}  Under the projection map $$(z,w_1,w_2,...,w_n,
v_1,v_2,...,v_n)\rightarrow(z,w_1,w_2,...,w_n)$$ complex tangents are
carried to complex tangents, so a complex tangent to $L^t$ at
$(z,w_1,w_2,...,w_n, v_1,v_2,...,v_n)$ must be of the form
$$(0,{\Cal T},{\Cal U})\tag 9$$ where ${\Cal T}$ is a complex tangent
to $K^t_z$ at $(w_1,w_2,...,w_n)$.  Then for (9) to be a complex
tangent, it will suffice that ${\partial I_j\over \partial \overline
{\Cal T}}=0$ for $j=1,2,...,n$.  Suppose, indeed, that ${\Cal
T}=(a_1,a_2,...,a_n)$ is a complex tangent and we map $\lambda
\mapsto (w_1,w_2,...,w_n)+\lambda(a_1,a_2,...,a_n)$.  Then
$\sum_{j=1}^na_j{\partial\rho\over \partial w_j}=0$, and for (9) to be
a complex tangent it will suffice that ${\partial I_j\over \partial
\overline\lambda}=0$ for all $j$.  We write $I(z,w)=(I_1(z,w),I_2(z,w),
...,I_n(z,w))={1\over S(z,w)}({\partial\rho\over
\partial w_1}(z,w),{\partial\rho\over
\partial w_2}(z,w),...,{\partial\rho\over\partial w_n}(z,w))$.
Then $$\eqalign{{\partial I\over\partial\overline \lambda}&=
{\partial({1\over S})\over \partial\overline \lambda}({\partial\rho\over
\partial w_1}(z,w),{\partial\rho\over
\partial w_2}(z,w),...,{\partial\rho\over\partial w_n}(z,w))\cr
&\hskip -1 mm+{1\over S}(\sum_{j=1}^n\overline{a_j}
{\partial^2\rho\over\partial w_1\partial \overline {w_j}}(z,w),
\sum_{j=1}^n\overline{a_j}{\partial^2\rho\over\partial w_2\partial 
\overline {w_j}}(z,w),
...,\sum_{j=1}^n\overline{a_j}{\partial^2\rho\over\partial w_n\partial
\overline {w_j}}(z,w)).\cr}\tag 10$$
If we take the complex inner product of (10) with
$(\overline{a_1},\overline{a_2}, ...,\overline{a_n})$ (the latter on
the right), then the first term drops out since
$\sum_{j=1}^na_j{\partial\rho\over \partial w_j}=0$.  The second term
becomes $${1\over S}\sum_{i,j=1}^n{\partial^2 \rho\over
\partial z_i\partial\overline{z_j}}a_i\overline{a_j}$$
which (exists and) is nonzero because $S\neq 0$ and $K^t_z$ is
strictly pseudoconvex.  We conclude that $L^t$ indeed has no nonzero
complex tangents. $\square$

{\bf Lemma 2.}  {\it $\widehat{(L^t)}_z=L^t_z$ for $z\in\Gamma$.}

{\it Proof.}  It is easy to show that $\widehat {(L^t)}_z=\widehat
{(L^t_z)}$ so all we must do is show that $\widehat {(L^t_z)}=L^t_z$.
Suppose that $(w',v')\in\widehat {L^t_z}$.  Then $w'\in \widehat
{K^t_z}$ since projection is analytic, so $\rho (z,w')\leq 1$.  Now
the polynomial $w_1v_1+w_2v_2+...+w_nv_n-1$ is identically $0$ on
$L^t$ so is $0$ on $\widehat {L^t_z}$.  Thus $w'\neq 0$.  Let ${\Cal
P}(w_1,w_2,...,w_n, v_1,v_2,...,v_n)=(w_1,w_2,...,w_n)$ and ${\Cal
Q}(w_1,w_2,...,w_n,v_1,v_2,...,v_n)=(v_1,v_2, ...,v_n)$.  Now for
$|w|$ sufficiently small, $T_w\equiv\{(v_1,v_2,...,v_n)\bigl|
w_1v_1+w_2v_2+...+w_nv_n=1\}$ doesn't meet $\widehat {{\Cal
Q}(L^t_z)}$.  Let ${\Cal B}_z=\{(w_1,w_2,...,w_n)\bigl |T_w\cap
\widehat {{\Cal Q}(L^t_z)}=
\emptyset\}$.  Then ${\Cal B}_z$ is open.  We claim that ${\Cal B}_z
\cap {\Cal P}
(\widehat{L^t_z})=\emptyset$.  Suppose $w'\in{\Cal P}
(\widehat{L^t_z})$, and $(w',v')\in\widehat{L^t_z}$ projects to it.
Then $v'\in {\Cal Q}(\widehat{L^t_z})\subset \widehat{{\Cal
Q}(L^t_z)}$ and $v'\in T_{w'}$ so $w'\notin {\Cal B}_z$, so the claim
holds.  Consider the connected component of ${\Cal B}_z\cap\{w\bigl |
\rho(z,w)<1\}$ which contains $0$.  Construct a continuous path
$w(s)$, $0\leq s\leq 1$, from $0$ to a point in the boundary of this
component.  Then we claim that $T_{w(1)}$ meets $\widehat {{\Cal
Q}(L^t_z)}$ in only one point, which must be in ${\Cal Q}(L^t_z)$.  If
they met in a point in $\widehat {{\Cal Q}(L^t_z)}\setminus{\Cal
Q}(L^t_z)$ then as $s
\rightarrow 1$, $T_{w(s)}$ approaches $\widehat {{\Cal Q}(L^t_z)}
\setminus{\Cal Q}(L^t_z)$, so by the Oka-Weil Theorem, must
also approach a point of ${\Cal Q}(L^t_z)$.  (What one must do is
consider the reciprocals of the complex affine maps which vanish on
$T_{w(s)}$.  They are analytic in a neighborhood of $\widehat{{\Cal
Q}(L^t_z)}$, so satisfy the maximum modulus principle with respect to
${\Cal Q}(L^t_z)$.  Clearly then $T_{w(s)}$ cannot approach a point of
$\widehat {{\Cal Q}(L^t_z)}\setminus{\Cal Q}(L^t_z)$ without
approaching a point of ${\Cal Q}(L^t_z)$.)  Hence $T_{w(1)}$ contains
a point of ${\Cal Q}(L^t_z)$, say $v'$.  But by strict hypoconvexity
of $\widehat {K^t_z}$, there is only one point $w'$ in
$\widehat{K^t_z}$ such that $\sum_{j=1}^nw'_jv'_j=1$; we must have
$w(1)=w'$ and $w'\in K^t_z$.  Thus ${\Cal B}_z$ is the entire interior
of $\widehat {K^t_z}$.  We conclude that ${\Cal P}(\widehat{L^t_z})$
is only $K^t_z$.

Then $\widehat{L^t_z}$ is a compact set in $K^t_z\times\C^n$.  Let
$(w',v')$ be the point in $L^t_z$ with first coordinate $w'$.  Suppose
that $(w',v'')\in\widehat{L^t_z}\setminus L^t_z$.  Let $Q$ be a
polynomial in $v$ which vanishes at $v'$ but equals $1$ at $v''$.
Choose a sequence of complex $(n-1)-$dimensional hyperplanes external
to $K^t_z$ approaching $w'$.  Suppose they are defined by the
vanishing of the complex affine functions $Q_i(w)$.  Regarding
$1/Q_i(w)$ as a function of both $w$ and $v$, let $M_i$ be the maximum
of $1/Q_i$ on $L^t_z$.  Then for large $i$, the modulus of
$Q(v){1\over (M_iQ_i(w))^i}$ is larger at $(w',v'')$ than at any
point on $L^t_z$, a contradiction since $Q(v){1\over (M_iQ_i(w))^i}$
is analytic on a large compact subset of $L^t_z\times\C^n$, hence is
approximable by polynomials.  This proves $\widehat{L^t_z}=L^t_z$.
$\square$

From Theorem 1, we have a parametrization of analytic graphs
with boundaries in the various $K^t$ by their values
at zero.  We wish to show that this parametrization is essentially
unique.  To do this, we need the following lemma.

{\bf Lemma 3. }  {\it Suppose that $w_0\in\C^n$, $N(w_0)$ an
open neighborhood such that there exist $F(w),G(w)$ with 
the properties arising from Theorem 1.  Then for a smaller
neighborhood $N'$ of $w_0$ the mapping
$P :\C^n\rightarrow\R$ given by $P (w)=\rho(z,F(w)(z))$
(where $z$ is any point in $\Gamma $)
is $C^1$ in $N'(w_0)$ with nonzero gradient.}

{\it Proof.}  Smoothness is trivial: $P $ is a composition of $C^1$
functions.  It will suffice to restrict ourselves to the case where
$g=(1,0,0,0,...,0)$ in Theorem 1, since, using notation from Theorem
1, $P(w)=\tilde P(w\cdot M(0))$, where $\tilde P(w)\equiv
\tilde\rho(z,\tilde F(w)(z))$ and $M(0)$ is invertible.  Let
$P'(u,v,k)(z)=\rho(z,F(u,v,k)(z))$, in the notation of Theorem 1.  Now
$DP'(0,0,0)[u,v,k]$ changes only with changes in $u$, and is injective
in $u$ (the proof is similar to an argument in the proof in Theorem 1
and the argument in section 2 of [F1]; $DP'(0,0,0)[u,v,k]=Iu$, where
$I(z)=2\Re(X(z){\partial\rho\over\partial w_1}(z,f(z)))$.)  Thus it
suffices to show that the derivative of $u(w)$ in $w$ is not
degenerate at $w_0$, and for this it suffices to show that the
derivative of the mapping $e(w)=u(w)(0)$ at $w_0$ is nondegenerate.
Begin with the equation
$w=f(0)+(u(w)(0))X(0)f(0)+v(w)(0)iX(0)f(0)+k(w)(0)$; taking the real
inner product of both sides with $X(0)\overline g(0)$, we obtain
$\Re\langle w,X(0)\overline g(0)\rangle =\Re
X(0)+\Re(u(w)(0))|X(0)|^2=\Re X(0)+e(w)|X(0)|^2$.  Differentiating
with respect to $w$, we find $\Re\langle w,X(0)\overline g(0)\rangle
=De(w_0)[w]|X(0)|^2.$ Letting $w=X(0)f(0)$ we find that
$De(w_0)[X(0)f(0)]=1\neq0$.  (This uses the fact that $X$ has winding
number zero, so $X(0)\neq 0$.)  Thus $De(w_0)$ is non-degenerate, as
desired.  $\square$

We henceforth assume that $N(w_0)$ has been shrunk so that it
possesses the property found in Lemma 3.

{\bf Lemma 4.}  {\it Suppose that $N(w_0),F(w),G(w)$ arise out of
Theorem 1 and Lemma 3 and $\phi $ is a function in $H^\infty(\Delta
)^n$ such that $\phi(z)=w$ for some $w\in N(w_0)$ and
$\rho(z,\phi(z))\leq P(w)$ for almost every $z\in\Gamma $.  Then
$F(w)\equiv\phi$.}

{\it Proof.}  If we have $F(w)\not\equiv\phi$, then the function
$\sum_{j=1}^n G_j(w)(F_j(w)-\phi_j)$ is not identically zero, because
$\sum_{j=1}^nG_j(w)(z)(F_j(w)(z)-\phi_j(z))\hskip -1mm=\hskip -1mm 
1-\sum_{j=1}^n
G_j(w)(z)\phi_j(z)$ which, for a.{}e.{} $z\in\Gamma $, can only equal
zero at $z\in\Gamma$ if $F(w)(z)=\phi(z)$ by the strict hypoconvexity
of $\rho $ and by the fact that $\rho(z,\phi(z))\leq P(w)$ for
a.{}e.{} $z\in\Gamma $.  However, if $F(w)\not\equiv \phi$ then
$F(w)(z)\neq \phi(z)$ on a set of positive measure in $\Gamma $.  Let
$t_0=P(w)\equiv\rho(z,F(w)(z))$.  By Lemma 3, we can choose a sequence
$\{w^j\}$, such that $w^j\rightarrow w$, $P(w^j)=t_j$ and
$t_j\downarrow t_0$.  Consider the function $z\mapsto 1-\sum_{j=1}^n
G_j(w^i)(z)\phi_j(z)$ for $i\geq 1$; we show in the next
paragraph that it has no zeroes on the disk (deforming
$G(w^i)$ through $g^t$ to a function near $0$ so that $\{v\bigl|
\sum_{j=1}^n g^t(z)v_j=1\}$ doesn't meet $\widehat {K^{t_0}_z}$ for
$z\in\Gamma$).  However, as $i\rightarrow\infty$,
$1-\sum_{j=1}^n G_j(w^i)\phi_j$ converges uniformly on compact sets to
$1-\sum_{j=1}^n G_j(w)\phi_j$, since $G$ is continuous in $w$.  By
Hurwitz' theorem this means that $1-\sum_{j=1}^nG_j(w)\phi_j$ is
identically zero on the disk since
$1-\sum_{j=1}^nG_j(w)(0)\phi(0)=1-\sum_{j=1}^nG_j(w)(0)F_j(w)(0)=0$.
(This holds because $F(w)(0)=w$ by Theorem 1.)  Thus we conclude that
$\sum_{j=1}^nG_j(w)(F_j(w)-\phi_j)$ is identically zero on $\Delta $,
so over $\Gamma $ in particular.  By the strict hypoconvexity of the
fibers, we conclude that $F(w)=\phi$ identically, as desired.

Now as to the deformation of $G(w^i)$ described, assuming $i$ fixed:
we can choose a homotopy $\{f^t\}$ of $F(w^i)$ such that
$\rho(z,f^t(z))\equiv t$ for all $t_i\leq t\leq R$ and $z\in\Gamma$,
$|f^t(z)|=t$ for $t\geq R$ and $z\in\Gamma $, and $f^{t_i}=F(w^i)$.
Then let $g^t(z)=D_w\rho(z,f^t(z))/\sum_j
f^t_j(z){\partial\rho\over\partial w_j} (z,f^t(z))$ for $t_i\leq t\leq
R$ and $\overline{f^t(z)}/t^2$ if $t>R$.  Extend $f^t,g^t$
harmonically to the closed disk.  If $R'$ is chosen large enough then
$|\sum_{j=1}^ng^{R'}_j(z)\phi_j(z)|\leq {1\over 2}$ for
$z\in\Delta$.  We then claim that for $r$ near $1$, there exists an
$\epsilon>0$ such that $|1-\sum_{j=1}^ng^t_j(z)\phi_j(z)|\geq
\epsilon$ for $t_i\leq t\leq R'$ and $r\leq |z|<1$.  If not, there
exists a convergent sequence $\{z_k,s_k\}$ such that $|z_k|\rightarrow
1$, $\phi(z_k)\rightarrow v$ and
$|1-\sum_{j=1}^ng^{s_k}_j(z_k)\phi_j(z_k)|
\rightarrow 0$.  Suppose that $(z_k,s_k)\rightarrow(z,s)$.
Then $g^{s_k}(z_k)\rightarrow g^s(z)$ and $v\in \widehat K^{t_0}_z$
(since the graph of $\phi $ is in the polynomial hull of $K^{t_0}$),
so $1-\sum_{j=1}^ng^{s}_j(z)v_j=0$. This is a contradiction because
$s>t_0$ and the set $\{w\in \C^n\bigl |
1-\sum_{j=1}^ng^{s}_j(z)w_j=0\}$ does not meet $\widehat{K^{t_0}_z}$.
With the existence of $r$ as claimed, we see that $1-\sum_{j=1}^n
G_j(w^i)\phi_j= 1-\sum_{j=1}^ng^{t_i}_j(z)\phi_j$ has the same winding
number on radius $r$ as $1-\sum_{j=1}^ng^{R'}_j(z)\phi_j$ which is
within ${1\over 2}$ of $1$ on $\Delta $.  Thus $1-\sum_{j=1}^n
G_j(w^i)\phi_j$ has no zeroes inside the disk of radius $r$, and none
outside by definition of $r$.  $\square$

{\bf Corollary.}  {\it If $N(w_0),F^0(w),G^0(w)$,
$N(w_1),F^1(w),G^1(w)$ both arise out of Theorem 1 and Lemma 3
such that $N(w_0)\cap N(w_1)\neq \emptyset$ then for all $w\in
N(w_0)\cap N(w_1)$, $F^0(w)=F^1(w)$ and $G^0(w)=G^1(w)$.}

{\it Proof.}  Choose $w\in N(w_0)\cap N(w_1)$.  Assume without loss of
generality that $\rho(z,F^1(w)(z))\leq \rho(z,F^0(w)(z))$.  Applying
Lemma 4, letting $\phi=F^1(w)$ and $F=F^0$, we have $F^0(w)=F^1(w)$,
from which  $G^0(w)=G^1(w)$ follows immediately. $\square$

We now explain in rough terms, without stating a precise theorem, how
we shall begin with an $f,g$ as in Theorem 1 and construct $F(w),G(w)$
for many $w$.  In our applications, the set where $\rho=R$ will be
equal to the set where $|w|$ is some constant $R$, so we could begin
with $f(z)=(R,0,0,...,0)$ and $g(z)=({1\over R},0,0,...,0)$.  Using
Theorem 1, we can construct a neighborhood $N$ of $f(0)$ and
associated $F(w),G(w)$.  Suppose $v$ is a boundary point of $N(f(0))$
and suppose $\{w^j\}$ is a sequence in $N$, $0<P(w^j)\leq R$, such
that $w^j\rightarrow v$ and $P(w^j)\rightarrow t>0$, with associated
$(F(w^n),G(w^n))$.  By Corollaries 1 and 2 in section 2 of [HMa], a
local uniform limit $(\phi,\psi)$ of the $(F(w^j),G(w^j))$ has a graph
whose accumulation points lie in $L^t$ (see the beginning of \S 2 for
the definition of $L^t$).  Then the functions $\phi,\psi$ are
$C^4$-smooth on $\Delta $ from [\v Ci] and they satisfy the properties
that $f,g$ do in Theorem 1; hence by Theorem 1, they may be
parametrized locally smoothly by $w\mapsto(F'(w),G'(w))$ in some
$N'(v)$, regarding $F'(w)$ and $G'(w)$ as elements of
$H^{1,2}(\Delta)^n$.  By the Corollary, where $N'(v)$ meets $N(v)$,
$(F',G')=(F,G)$.  Thus we can extend $F,G$ as far as the graphs of the
limiting functions $(F(v),G(v))$ continue to have boundaries where
$\rho $ is $C^6$ and strictly hypoconvex, and $0<\rho\leq R$.

If $\rho $ has smoothness $C^k$ for $k>6$ then the various $F(w)$ and
$G(w)$ extend to $C^{k-2}(\Gamma)$, again using \v Cirka's result,
since then $L^t$ is a totally real $C^{k-1}$ manifold.
 
\S 3 {\bf Polynomial hulls with hypoconvex fibers.}
  
{\bf Theorem 2.}  {\it Suppose $\rho $ satisfies (3) with $S=0$, so
the fibers of $K^t_z$ of $K^t$ enclose the origin in $\C^n$.  Then
$\widehat K\setminus K$ is the union of graphs of elements of
$A(\Delta )\cap C^4(\Gamma)$ (whose boundaries lie in some $K^t$,
$t\leq 1$).  Given a point in $\partial \widehat
K\cap(\opdisk\times\C^n)$, there is precisely one element of
$H^\infty(\Delta )$ whose graph is in $\widehat K\setminus K$ and
passes through that point.  For all $z\in\Delta $, $\widehat K_z$ is
hypoconvex, with $C^1$ boundary; in fact, we have $\widehat
K_0=\{w\in\C^n\bigl| P(w)\leq1\}$.}

{\it Proof.}  We easily find $f,g$ as in Theorem 1; we can take
$f(z)=(R,0,0,...,0)$ and $g(z)=({1\over R},0,0,...,0)$.  Then we can
construct an open $U$ and $F,G,P$ as before.  We claim that we can
extend $F,G$ and $P$ smoothly to $\{w|0< |w|\leq R\}$ and $P$
continuously to where $w=0$.  The extension of $F,G$ to the set where
$|w|=R$ is obvious, and then the extension to a neighborhood of this
sphere is by application of Theorem 1 and Lemmas 1-4.  Then $P(w)=R$
if $|w|=R$, so $D_w P(w)$ is a real multiple of $w$ (it is nonzero, by
Lemma 3); we claim it is a positive multiple.  Were it negative, then
for some $|w|>R$, we would have $|F(w)|\leq R$ on $\Gamma $, but
$|F(w)(0)|>R$, in violation of the maximum modulus principle.  Suppose
$U$ is the maximal open subset of $\{w|0< |w|<R\}$ to which $F$ and
$G$ (so $P$ also) can be extended $C^1-$smoothly.  If $U$ excludes
points in the annulus $0<|w|\leq R$, then take an open segment in $U$
one of whose vertices $w$ is not in $U$, $0<|w|<R$ and the other of
which lies on $\{w||w|=R\}$.  We claim that $P$ is bounded on the
segment.  If not, then at some point $w'$ on the segment $P(w')=R$, so
$F(w')$ and $G(w')$ are constant functions of modulus $R,{1\over R}$,
respectively, so $|w|=R$, a contradiction.  (Recall that $\overline
{G(w')(z)}$ is a complex multiple of a normal to the sphere of radius
$R$, in fact, it equals ${1\over R^2}F(w')(z)$; then $1=\langle
F(w'),\overline {G(w')}\rangle ={1\over R^2} |F(w')|^2$ is bounded
analytic and real on $\Gamma $.  Thus $F(w')$ has constant modulus on
$\Delta $, so its components are constants, from which we can make the
above conclusions.)  Then take a sequence $\{w^j\}$ on the segment
converging to $w$, such that $\{P(w^j)\}$ converges.  We cannot have
$\lim _{j\rightarrow\infty}P(w^j)=0$ because then a subsequence of
$\{F(w_j)\}$ converges locally uniformly to zero, which is impossible
because $w^j=F(w^j)(0)$ converges to $w$ which is not zero.  Then we
can define $F(w),G(w)$ as the local uniform limits of subsequences of
$\{F(w_j)\}$ and $\{G(w_j)\}$.  As outlined at the end of \S 2, $F(w)$
and $G(w)$ are $C^4$ functions on $\Gamma $.  Then using Theorem 1, we
extend smoothly to a neighborhood of $w$ in such a way as to coincide
with $F$ and $G$ on $U$, a contradiction.  Thus $F,G$ extend to be
$C^1$ on the annulus $0< |w|\leq R$, where $F,G$ must be the natural
constants when $|w|=R$: $F(w)(z)=w$ constantly and $G(w)(z)=\overline
w/|w|^2$.

Suppose $\{w^j\}$ is a sequence which converges to the origin
${\buildrel \rightarrow\over 0}$ but $t\equiv\lim P(w^j)\neq 0$.  Then
since $t>0$, we could define $F({\buildrel \rightarrow\over 0})$ and
$G({\buildrel \rightarrow\over 0})$ to be local uniform limits of
$\{F(w^j)\}$ and $\{G(w^j)\}$.  However, since $F({\buildrel
\rightarrow\over 0})(0)={\buildrel \rightarrow\over 0}$ we cannot
possibly have $\sum_{j=1}^nF_j({\buildrel \rightarrow\over 0})(0)
G_j({\buildrel
\rightarrow\over 0})(0)$ equal to $1$.  Thus $\lim P(w^j)= 0,$ which
means we can extend $P$ continuously to $\{w\bigl| |w|\leq R\}$ by
defining $P({\buildrel \rightarrow\over 0})=0$.  We also find from
this limit that as $w^j\rightarrow 0$, $F(w^j)$ converges uniformly to
zero, so for convenience we define $F({\buildrel \rightarrow\over
0})={\buildrel \rightarrow\over 0}$.

Let $s$ be the maximum of $P$ on $\widehat K_0$.  (Note that the
domain of $P$ clearly contains $\widehat
K_0\subset\{w\in\C^n\,\bigl|\,|w|\leq R\}$.)  Clearly $s<R$.  Suppose
$s>1$.  Let $v$ be a point on $\widehat K_0$ where this maximum is
attained.  Then $v$ is not in the interior of $\widehat K_0$ because
$P$ does not attain local maxima where it is smooth, since $P$ has
nonzero gradient.  The only place $P$ is not smooth when $|w|<R$
occurs when $w=0$ and $P$ doesn't attain a local maximum there since
$P({\buildrel \rightarrow\over 0})=0$.  Thus we can choose a
continuous path $v(t)$, $s<t<R$, outside of $\widehat K_0$ such that
$P(v(t))=t$ and as $t\rightarrow s^+$, $v(t)\rightarrow v$.

Then consider 

$$\{(z,w)\in\Delta\times\C^n\bigl| \sum_{j=1}^n G_j(v(t))w_j=1\}
\tag 11$$ Let $s'$= the supremum of all $t$ such that (11) meets
$\widehat K$.  Clearly $1<s\leq s'<R$.

The function
$$M_t(z,w)={1\over \sum_{j=1}^n G_j(v(t))w_j-1}$$ on $\Delta
\times \C$ is defined on $\widehat K$ for $t>s'$.  Since
$G(v(t))\in A(\Delta)$ for all $t$, $M_t(z,w)$ is uniformly
approximable on $\widehat K$ for $t>s'$ by functions analytic in a
neighborhood of $\widehat K$, hence uniformly approximable by
polynomials in a neighborhood of $\widehat K$ by the Oka-Weil
Theorem. This means that
$$\sup _{(z,w)\in \widehat K}|M_t(z,w)|\leq\sup_{(z,w)\in
K}|M_t(z,w)|,$$ for $t>s'$.  As $t\downarrow s'$, $\sup _{(z,w)\in
\widehat K}|M_t(z,w)|\rightarrow \infty$ by the definition of $s'$.  
However since $s'>1$, the
distance between $K^{t}$ and $K^1$ is bounded away from zero uniformly
in $t$, and the singularity set of $M_t$ is no closer to $K^1=K$ than
points in $K^t$, by strict hypoconvexity of the fibers, so
$$\sup_{(z,w)\in K}|M_t(z,w)|$$ is bounded uniformly in $t$.  This
contradicts the previous assertion.  Thus $s=1$ and we find that
$\widehat K_0\subset \{w\in \C^n \bigl| |w|<R, P(w)\leq 1\}$.  Thus
every point in $\widehat K_0$ is on the graph of some $F(w)$ for which
$P(w)\leq 1$.  The same holds for all other $\widehat K_z$,
$z\in\opdisk$ by applying a M\" obius transformation to the disk
sending $z$ to $0$ and applying the same argument.  Thus $\widehat K$
is indeed the union of graphs of elements of $A(\Delta )$ which extend
to $C^4(\Gamma )$.

Now given a point in $\partial\widehat K\cap(\opdisk\times\C^n)$,
suppose by applying a M\" obius transformation that it has the form
$(0,w)$.  Then $P(w)=1$, for if $P(w)<1$, the graphs of $F(w)+e$ are
in $\widehat K$ for constant analytic vector valued functions $e$ of
sufficiently small modulus, so $(0,F(w)(0))$ is in the interior of
$\widehat K$.  Since $P(w)=1$, applying Lemma 4, there is no $\phi\in
H^\infty(\Delta )^n$ other than $F(w)$ such that $\phi(0)=w$ whose
graph is contained in $\widehat K\setminus K$, as $\phi(z)\in
\widehat K_z$ for almost every $z\in\Gamma $ for such a $\phi$.

The last statement of the theorem is already known for $z\in\Gamma $.
If $z\in\opdisk$, we may assume without loss of generality that $z=0$
by applying a M\" obius transformation.  If $P(w')\leq 1$, then
$w'\in\widehat K_0$; the converse was shown above.  Thus
$\partial\widehat K_0=
\{w\in\C^n\bigl | P(w)=1\}$ and $\widehat K_0$ has $C^1$
boundary since $P$ is $C^1$ for $0<|w|\leq R$.  To show that $\widehat
K_0$ is hypoconvex: suppose that point $w'\notin
\widehat K_0$.  If $|w'|\leq R$ then the complex affine hyperspace
$\{w\in\C^n\,|\,\langle G(w')(0),\overline w\rangle =1\}$ is external
to $\widehat K_0$ (since $P(w')>1$) and passes through $w'$.  If
$|w'|>R$, it is easy to find such a hyperspace.  Thus $\widehat K_0$
is hypoconvex.  $\square $

Theorem 2 will hold if in (3), $S$ has polynomial coordinates, so the
sets $K^t_z$, instead of enclosing the origin, encircle points $S(z)$.

\S 4 {\bf The $H^\infty$ control problem.}

Assume that $\rho $ satisfies (3), let $\gamma_\rho$ be as defined
in (2) and let $\delta_\rho$ be defined by
$$\delta_\rho\equiv\inf_{f\in A(\Delta)^n}\sup_{z\in\Gamma}
\rho(z,f(z)).$$  We assume that $S$ is not analytic, so that 
$0<\gamma_\rho\leq\delta_\rho$.  For $m=1,2,3,...,$ let
$\rho_m(z,w_1,w_2,...,w_n)=
\rho(z,{w_1\over z^m},{w_2\over z^m},...,{w_n\over z^m})$ and define
$z^mK^t$ to be the set in $\Gamma\times\C^n$ where $\rho_m$ equals
$t$.  Now since $\rho_m(z,z^mg(z))=\rho(z,g(z))$ for all $g\in
H^\infty(\Delta )^n$, we find that $\delta_{\rho_m}\leq\delta_\rho$
and $\gamma_{\rho_m}\leq\gamma_\rho.$ We claim that there exists an
$m>0$ such that $\delta_{\rho_m}<\delta_\rho$.  The reason for this is
as follows: choose any continuous selector $\alpha $ for
$\{(z,w)\in\Gamma\times\C^n\bigl|
\rho(z,w)\leq
\delta_\rho-\epsilon\}$, where of course $\epsilon $ is chosen small
enough so that this is possible.  Then there exists an $n$-tuple of
harmonic polynomials $q=(q_1,q_2,...,q_n)$ such that $q$ is so close
to $\alpha $ uniformly, that $q$ is a harmonic polynomial selector for
$\{(z,w)\in\Gamma\times\C^n\bigl|
\rho(z,w)\leq
\delta_\rho-{\epsilon\over 2}\}$.  There exists an $m>0$ such that
$z^mq(z)$ has coordinates which are analytic polynomials.  Hence
$\delta_{\rho_m}\leq\delta_\rho-{\epsilon\over 2}$.  Let $m$ in fact
be the least positive integer such that $\delta_{\rho_m}<\delta_\rho$.
Then there exists an $n$-tuple of harmonic polynomials $q(z)$ such
that $\rho(z,q(z))<\delta_\rho$ and $z^mq(z)$ is analytic.  Note that
this means $m\geq 1$.

We shall not assume that the level sets of $\rho $
enclose the origin as we did earlier.  The role of the zero
graph is here replaced by the graph of $q$.

{\bf Theorem 3.}  {\it Suppose that $\rho $ satisfies (3) with
$S(z)=q(z)$.  Then there exists a unique
$\phi\in H^\infty(\Delta )^n$ for which $\rho(z,\phi(z))\leq
\gamma_\rho$, for which in fact $\phi$ extends to be $C^4$ on $\Gamma
$ and $\rho(z,\phi(z))=\gamma_\rho$ for all $z\in\Gamma$.  Also, if
$\rho(\overline z,\overline w)=\rho(z,w)$, then $\phi $ is
\R$^n$-valued on the real axis, i.e. $\overline {\phi(\overline
z)}=\phi(z)$.  If $\rho $ is $C^k$ where $\rho\neq 0$ for $k>6$
then $\phi $ extends to $C^{k-2}(\Gamma)$.}

The condition of being real on the real axis has applications in
engineering.  See [HMe2].  The condition that
$\rho(z,\phi(z))=\gamma_\rho$ for all $z$ is known as the ``frequency
domain Bang-Bang principle,'' or a condition of ``flat performance''.
The existence of an $H^\infty$ solution in the theorem has already
been proven in [HMa].  If the $H^\infty$ solution is smooth on $\Gamma
$, then it is the only smooth solution; see [V].

{\it Proof.}  From the definition of $\rho $ it is clear that
$\gamma_\rho\leq R$ and $\delta_\rho\leq R$.

Let us suppose that $m=1$, i.e.,
$\delta_{\rho_1}<\delta_\rho$.  
Applying the technique of Theorem 2 to $
\rho_1(w)$ (replacing the zero function by $p(z)=zq(z)$), 
we obtain $F(w)$ and $P(w)$ to be defined and $C^1$ in the region
where $|w|\leq R+\epsilon$ except when $w=p(0)$.  Note that $p(0)\neq
{\buildrel \rightarrow\over 0}$ since $q$ is not analytic.  Thus $P$
and $F$ are defined at the origin.

We claim that $\phi\equiv{1\over z}F({\buildrel \rightarrow\over 0})$
will solve the $H^\infty$ control problem (2) for $\rho $ and
$\gamma_\rho=P({\buildrel \rightarrow\over 0})$.  Clearly $\phi\in
A(\Delta)\cap C^4(\Gamma )$ because $F({\buildrel \rightarrow\over
0})={\buildrel \rightarrow\over 0}$ and $F({\buildrel \rightarrow\over
0})\in A(\Delta)\cap C^4(\Gamma )$.  Also, if $\rho $ is in fact $C^k$
for $k>6$ then from the observation at the end of \S 2, we also find
that $\phi\in C^{k-2}(\Gamma)$.  The boundary values of $\phi $ are
clearly in $K^{P({\buildrel \rightarrow\over 0})}$ since the same
holds for $F({\buildrel \rightarrow\over 0})$ with respect to
$zK^{P({\buildrel
\rightarrow\over 0})}$.  We can show that there is no other element
$\psi $ in $H^\infty(\Delta )$ for which
ess$\,\sup_{z\in\Gamma}\rho(z,\psi(z))\leq P({\buildrel
\rightarrow\over 0})$ because then there would exist another element
$z\psi $ in $H^\infty(\Delta )$ with value ${\buildrel
\rightarrow\over 0}$ at $z=0$ for which $\rho_1(z,z\psi(z))\leq
\gamma_\rho $ for almost every $z\in\Gamma $. This contradicts Lemma 4.

To prove that $\phi(z)=\overline{\phi(\overline z)}$ if $\rho(z,w)=\rho(
\overline z,\overline w)$, we imitate [HMa]: note that
$\overline{\phi(\overline z)}$ is analytic in $z$ and
$\rho(z,\overline{\phi(\overline z)})=\rho(\overline z,\phi(\overline
z)) =\gamma_\rho$ for all $z\in\Gamma $.  Thus
$\overline{\phi(\overline z)}$ is another solution to the $H^\infty$
control problem, so must be the same as $\phi(z)$ by the uniqueness
just proven.

Now let us assume that $m>1$.  Applying the work above to
$\rho_{m-1}$, we find there is a unique solution $k$ to the $H^\infty
$ control problem for $\rho_{m-1}$ which also happens to be in
$A(\Delta )\cap C^4(\Gamma)$.  Now $\gamma_\rho \leq\delta_\rho=
\delta_{\rho_{m-1}}=\gamma_{\rho_{m-1}}$.  (The first equality
is by definition of $m$; the second because we showed the $H^\infty$
solution for $\rho_{m-1}$ to be in $A(\Delta )$.)  However, we already
know $\gamma_{\rho_{m-1}}\leq \gamma_\rho,$ so
$\gamma_{\rho_{m-1}}=\gamma_\rho.$ Thus any solution $f$ to the
$H^\infty $ control problem for $\rho $ must satisfy
$z^{m-1}f(z)=k(z)$ for all $z\in\opdisk$.  Using the known existence
of a solution to (2) from [HMa], this proves uniqueness, smoothness of
the solution and flatness of performance for the $H^\infty $ control
problem for $\rho $.  The property of being
\R$^n$-valued on the real axis follows as well since $z^n$ is real on
the real axis.  $\square $

Using a version of Theorem 1, it is easy for us to show that if $\rho$
varies smoothly then the solution to the $H^\infty $ control problem
and the optimal control also vary smoothly.

We suppose that $\rho:\Gamma\times\C^n\times I$ and $S(z,\tau)$ are
$C^6$ where $I$ is an open interval in $\R$ and for every $\tau\in I$,
$\rho^\tau(z,w)\equiv\rho(z,w,\tau)$ satisfies (3) with respect to
$S^\tau(z)\equiv S(z,\tau)$.

{\bf Theorem 4.}  {\it Suppose that there exists an $m\geq 1$ which
for all $\tau\in I$ is the least positive integer such that
$\delta_{\rho^\tau_m}<\delta_{\rho^\tau}$ and suppose that for all
$\tau\in I$, $z^mS^\tau(z)$ is in $A(\Delta)$ but $z^{m-1}S^\tau(z)$
is not.  If $H(\tau)$ denotes the solution to the $H^\infty$ control
problem for $\rho^\tau,$ then $H:I\rightarrow H^{1,2}(\Delta )$ and
$\gamma_{\rho^\tau}$ are $C^1$ functions of $\tau$.}

{\it Proof.}  Fix any point in $I$, say $0$ without loss of
generality.  We consider the function
$\rho_m(z,w,\tau)\equiv\rho(z,{w\over z^m},
\tau)$ on $\Gamma\times\C^n\times I$.  We can use 
reasoning similar to Theorem 1 and Lemmas 1-4, to conclude the
existence of $C^1$ function $F(w,\tau)$ in $N(0)\times I$ such that,
following Theorem 3, ${1\over z^m}F({\buildrel \rightarrow\over
0},\tau)$ is the solution to the $H^\infty$ control problem for
$\rho^\tau$.  (The smoothness of the associated $\Phi$ follows from
the Lemma in the Appendix, replacing $n$ by $n+1$ and regarding $\tau
$ as a constant function $f_{n+1}$.)  Then it is easy to see that
$H(\tau)={1\over z^m}F({\buildrel\rightarrow \over 0},\tau)$ and
$\gamma_{\rho^\tau}=
\rho(1,H(\tau)(1),\tau)$ are $C^1$ functions of $\tau $. $\square $

\S 5 {\bf Appendix.}

As promised, we prove that $\Phi $ is a $C^1$-differentiable map in a
neighborhood of $(0,0,0,0,w_0)$.  Since harmonic conjugation is
continuous linear on $W^{1,2}_R(\Gamma )$ (so smooth), $\tilde u\in
W^{1,2}_R(\Gamma )$ if $u$ is.  The product and chain rule for Sobolev
functions show that $\Phi_1(u,v,k,l,w)$ is in $W^{1,2}(\Gamma)^n$. We
may similarly conclude that $\Phi_2(u,v,k,l,w)\in W^{1,2}(\Gamma) ^n$,
where the only additional facts that we need is that the denominator
is bounded away from $0$ for small $u,k$; this holds because if $u,k$
are small in $W^{1,2}$ then they are uniformly near $0$, so that
denominator is uniformly near $\sum_{j=1}^n
\overline {f_j(z)}{\partial\rho\over \partial \overline {w_j}}(z,f(z))$,
which is never zero for $z\in\Gamma$.  Then using the fact that
$W^{1,2} $ is an algebra, we conclude that indeed
$\Phi_2(u,v,k,l,w)\in W^{1,2}(\Gamma) ^n$.  If the integrand
of $\Phi_4$ is $C^1$, so is $\Phi_4$ since the integral is
continuous linear.

It will then suffice to prove the following lemma.

{\bf Lemma. }  {\it Let $p$ be a $C^3$ function on $\Gamma\times\C^n$.
If $P:(W^{1,2}(\Gamma))^n\rightarrow W^{1,2}(\Gamma)$ is given by
$P(f)(z)\equiv p(z,f(z))$, then $P$ is a $C^1$ function.}

{\it Proof. }  We claim that $DP(f)$ is given by the map $T_f\in L(
(W^{1,2}_C)^n,W^{1,2}_C)$ such that $T_f[h](z)=Dp(z,f(z))[0,h(z)]$.
We must first check that $T_f$ is in the desired space.  This is not
difficult: if $h=(h_1,h_2,...,h_n)$, then $T_f[h]=\sum_{j=1}^n
r_jh_j+s_j\overline{h_j},$ where $r_j(z)={\partial p\over\partial
w_j}(z,f(z)), s_j(z)={\partial p\over\partial \overline{w_j}}(z,f(z))$
are $W^{1,2}$ functions.  It is then a simple exercise to show that
$\|T_f[h]\|_{1,2}$ is less than or equal to a constant times
$\|h\|_{1,2}$.  It is also a simple matter to show that $T_f$ varies
continuously in $f$.  Next, by Taylor's formula,
$P(f+h)(z)-P(f)(z)=p(z,f(z)+h(z))-p(z,f(z))=Dp(z,f(z))[0,h(z)]
+R(f,h)(z),$ where
$$R(f,h)(e^{i\theta})=\int_0^1(1-t)D^2p(e^{i\theta},f(e^{i\theta})+
th(e^{i\theta}))[0,h(e^{i\theta})][0,h(e^{i\theta})]\,dt.\tag 12$$ For
$h$ with small norm in $W^{1,2}$, $h$
also has small supremum norm $\leq C\|h\|_{1,2}$; Then
$D^2p(e^{i\theta},f(e^{i\theta})+ th(e^{i\theta}))$ is uniformly
bounded in $\theta,t$ for such $h$, so the above integral has absolute
value bounded by a constant times $\|h\|_\infty^2$.  Thus
$\|R(f,h)\|_2\leq C(f)\|h\|_\infty^2\leq C_1(f)\|h\|^2_{1,2},$ so
$$\lim_{h\rightarrow 0\,in\, W^{1,2}}{\|R(f,h)\|_2\over
\|h\|_{1,2}}=0.\tag 13$$  We may
differentiate under the integral sign in (12) to get an integrand so
that $|R(f,h)_\theta(e^{i\theta})|$ is bounded above by a constant
times
$|f_\theta|\|h\|_\infty^2+\|h\|_\infty^2|h_\theta|+\|h\|_\infty|h_\theta|$ 
$=\|h\|_\infty
(|f_\theta|\|h\|_\infty+\|h\|_\infty|h_\theta|+|h_\theta|)$ which is
less than or equal to a constant times
$\|h\|_{1,2}(|f_\theta|\|h\|_\infty+\|h\|_\infty|h_\theta|+|h_\theta|)$.
Then $\|R(f,h)_\theta\|_2\leq C_2\|h\|_{1,2}(\|
|f_\theta|\|h\|_\infty+\|h\|_\infty|h_\theta|+|h_\theta| \|_2)$
and we conclude that
$$ \lim_{h\rightarrow 0\,in\, W^{1,2}}{\|R(f,h)_\theta\|_2\over
\|h\|_{1,2}}=0,$$
since $f_\theta$ is in $L^2_C(\Gamma)$.  Combining this with (13), we
conclude that $P$ is differentiable.  We showed above that $DP(f)$ is
continuous in $f$, so we are done.  $\square$

{\bf Acknowledgments}

We would like to extend thanks to Professors Harold Boas, Al Boggess,
Emil Straube and John Wermer.

\enddocument